\documentclass[notitlepage,leqno,11pt]{amsart}
\usepackage{amsmath,amssymb,amsbsy,amsfonts,amsthm,latexsym,
            amsopn,amstext,amsxtra,euscript,amscd}
\usepackage{hyperref}

\begin{document}

\newtheorem{theorem}{Theorem}
\newtheorem{lemma}[theorem]{Lemma}
\theoremstyle{remark}
\newtheorem{remark}[theorem]{\bf Remark}

\title{A variation on the theme of Nicomachus}
%
\author{Florian Luca}
\address{School of Mathematics, University of the Witwatersrand, Private Bag X3, Wits 2050, Johannesburg, South Africa}
\address{Max Planck Institute for Mathematics, Vivatsgasse 7, 53111 Bonn, Germany}
\address{Department of Mathematics, Faculty of Sciences, University of Ostrava,
30 dubna 22, 701 03 Ostrava 1, Czech Republic}
\email{Florian.Luca@wits.ac.za}

\author{Gerem\'\i as Polanco}
\address{School of Natural Science, Hampshire College,
893 West St, Amherst, MA 01002, USA}
\email{gpeNS@hampshire.edu}

\author{Wadim Zudilin}
\address{IMAPP, Radboud Universiteit, PO Box 9010, 6500 GL Nijmegen, The Netherlands}
\email{w.zudilin@math.ru.nl}
\address{School of Mathematical and Physical Sciences, The University of Newcastle, Callaghan, NSW 2308, Australia}
\email{wadim.zudilin@newcastle.edu.au}

\begin{abstract}
In this paper, we prove some conjectures of K.~Stolarsky concerning the first and third moment of the Beatty sequences with the golden section and its square.
\end{abstract}

\maketitle

\section{Introduction}
\label{sec1}

Nicomachus's theorem asserts that the sum of the first $m$ cubes is the square of the $m$th triangular number,
\begin{equation}
\label{eq:nico}
1^3+2^3+\dots+m^3=(1+2+\dots+m)^2.
\end{equation}
In the notation
$$
Q(\alpha,m):= \frac{\sum_{n=1}^m \lfloor\alpha n\rfloor^3}{\bigl(\sum_{n=1}^m \lfloor\alpha n\rfloor\bigr)^2},
$$
where $\alpha\in\mathbb R\setminus\{0\}$, it implies that
\begin{equation}
\lim_{m\to\infty}Q(\alpha,m)=\alpha.
\label{nico-lim}
\end{equation}
Here, $\lfloor x\rfloor$ is the integer part of the real number $x$. The limit in \eqref{nico-lim} follows from $\lfloor\alpha n\rfloor=\alpha n+O(1)$
and Nicomachus's theorem~\eqref{eq:nico}.

Recall that the Fibonacci and Lucas sequences $\{F_n\}_{n\ge 0}$ and $\{L_n\}_{n\ge 0}$ are given by $F_0=0$, $F_1=1$ and $L_0=2$, $L_1=1$ and the recurrence relations
\begin{equation*}
F_{n+2} = F_{n+1}+F_n, \qquad
L_{n+2} = L_{n+1}+L_n
\end{equation*}
for $n\ge 0$. In an unpublished note, which was recorded by the second author,
Stolarsky observed that the limit relation \eqref{nico-lim} can be `quantified' for $\alpha=\phi$ and $\phi^2$,
where $\phi:=(1+{\sqrt{5}})/2$ is the golden mean, and a specific choice of $m$ along the Fibonacci sequence.
The corresponding result is Theorem~\ref{conj1} below. We complement it by a general analysis of moments of the Beatty sequences
and give a solution to a related arithmetic question in Theorem~\ref{conj6}.

\section{Principal results}
\label{sec2}

\begin{theorem}
\label{conj1}
For $k\ge 1$ an integer, define $m_k:=F_k-1$.
We have
\begin{align*}
Q(\phi^2,m_{2k})-Q(\phi,m_{2k})
&=\begin{cases}
1-\dfrac{1}{(F_{k+1})^2 L_{k+2} L_{k-1}} &\text{if $k$ is even}, \\[5pt]
1-\dfrac{1}{(L_{k+1})^2 F_{k+2} F_{k-1}} &\text{if $k$ is odd},
\end{cases}
\\ \intertext{and}
Q(\phi^2,m_{2k-1})-Q(\phi,m_{2k-1})
&=\begin{cases}
1-\dfrac{F_{k-2}}{F_{k+1} (F_k)^2  (L_{k-1})^2} &\text{if $k$ is even}, \\[5pt]
1-\dfrac{L_{k-2}}{L_{k+1}(L_{k})^2 (F_{k-1})^2} &\text{if $k$ is odd}.
\end{cases}
\end{align*}
\end{theorem}

The theorem motivates our interest in the numerators and denominators of $Q(\phi,F_k-1)$ and $Q(\phi^2,F_k-1)$, which can be thought as
expressions of the form
$$
A(k,s):=\sum_{n=1}^{F_k-1}\lfloor\phi n\rfloor^s
\qquad\text{and}\qquad
A'(k,s):=\sum_{n=1}^{F_k-1}\lfloor\phi^2n\rfloor^s
$$
for $k=1,2,\dots$ and $s=1,3$. Our analysis in Section~\ref{sec3} covers more generally the sums
\begin{equation}
A(k,s,j):=\sum_{n=1}^{F_k-1} n^j \lfloor\phi n\rfloor^s, \qquad\text{where}\quad k=1,2,\dots \quad\text{and}\quad s,j=0,1,2,\dotsc.
\label{Aksj}
\end{equation}
Namely, we find a recurrence relation for $A(k,j,s)$ and deduce the recursions for $A(k,s)=A(k,s,0)$ and $A'(k,s)$ from it.
The strategy leads to the following expressions for the desired numerators and denominators in Theorem~\ref{conj1}, which are given in Lemmas~\ref{conj2}--\ref{conj4}.

\begin{lemma}
\label{conj2}
Let $k\ge 1$ be an integer. Then
\begin{equation}
\label{eq:conj2}
\begin{aligned}
A(k,1) & = \frac{1}{2}(F_{k+1}-1)(F_{k}-1), \\
A'(k,1) & = \frac{1}{2}(F_{k+2}-1)(F_{k}-1).
\end{aligned}
\end{equation}
\end{lemma}

\begin{lemma}
\label{conj3}
Let $k\ge 1$ be an integer. Then
\begin{align*}
A(2k,3) &=  \frac{1}{4}(F_{2k-1}-1)(F_{2k+1}-1)^2(F_{2k+2}-1)
\\ \intertext{and}
A(2k-1,3)&=\frac14(F_{2k-1}-1)(F_{2k}-1)\times\frac15(L_{4k}-3L_{2k+1}-L_{2k}+3).
\end{align*}
\end{lemma}

\begin{lemma}
\label{conj4}
Let $k\ge 1$ be an integer. Then
\begin{align*}
A'(2k,3)&=\frac14(F_{2k}-1)(F_{2k+2}-1)\times\frac15(L_{4k+4}-5L_{2k+3}+13)
\\ \intertext{and}
A'(2k-1,3)&=\frac14(F_{2k-1}-1)(F_{2k+1}-1)\times\frac15(L_{4k+2}-5L_{2k+2}+7).
\end{align*}
\end{lemma}

Finally, we present an arithmetic formula inspired by Stolarsky's original question.

\begin{theorem}
\label{conj6}
For $k\ge 1$, we have
$$
\operatorname{LCM} \bigl(A(2k,1),A'(2k,1)\bigr)
=\begin{cases}
\dfrac{1}{2} F_{k+1} F_k L_{k+2} L_{k+1}  L_{k-1} & \text{if}\; 2\mid k, \\[4.5pt]
\dfrac{1}{2} F_{k+2} F_{k+1} F_{k-1}  L_{k+1} L_k & \text{if}\; 2\nmid k.
\end{cases}
$$
\end{theorem}

\begin{remark}
Lemmas \ref{conj2}--\ref{conj4} indicate that the expression
$$
Q(\phi^2,F_k-1)-Q(\phi,F_k-1)
$$
is expressible as a fraction whose numerator and denominator are polynomials in Fibonacci and Lucas numbers with indices depending linearly on $k$ according to the parity of $k$,
yet the statement of Theorem~\ref{conj1} presents formulas for these quantities according to the congruence class of $k$ modulo $4$ rather than modulo~$2$.
The discrepancy is related to different factorizations of the factors that occur in the formulas for $A(k,j)$ and $A'(k,j)$ for $j\in\{1,3\}$,
as each of such factors $F_n-1$ is a product of a Fibonacci and Lucas number according to the congruence class of $n$ modulo~$4$ (see formulas~\eqref{eq:fact}).
\end{remark}

\section{Recurrence relations for auxiliary sums}
\label{sec3}

Here, we show how to compute the integer-part sums \eqref{Aksj}.
This clearly covers the cases $A(k,s)=A(k,s,0)$. On using
$$
\phi^2=1+\phi,
$$
which upon multiplication with the integer $n$ and taking integer parts becomes
$$
\lfloor\phi^2n\rfloor= n+\lfloor\phi n\rfloor,
$$
one also gets the explicit formulas
$$
A'(k,s)=\sum_{i=0}^s\binom siA(k,s-i,i).
$$

Using the Binet formula
$$
F_k=\frac{\phi^k-(-\phi^{-1})^{k}}{\sqrt{5}}\quad {\text{\rm for~all}}\quad k\ge 0,
$$
one proves easily that
\begin{equation*}
\lfloor\phi F_k\rfloor=F_{k+1}-\epsilon_k,
\qquad\text{where}\quad \epsilon_k=\frac{1+(-1)^k}2,
\end{equation*}
and that
\begin{equation*}
\lfloor\phi(F_k+n)\rfloor=F_{k+1}+\lfloor\phi n\rfloor
\qquad\text{for}\quad 1\le n\le F_{k-1}-1
\end{equation*}
(see, for example, \cite{Ki}).  Thus,
\begin{align*}
A(k+1,s,j)
&= \sum_{n=1}^{F_k-1} n^j\lfloor\phi n\rfloor^s+F_k^j\lfloor\phi F_k\rfloor^s+\sum_{n=F_k+1}^{F_{k+1}-1} n^j\lfloor\phi n\rfloor^s
\\
&= A(k,s,j)+F_k^j(F_{k+1}-\epsilon_k)^s+\sum_{n=1}^{F_{k+1}-F_k-1} (F_k+n)^j\lfloor\phi(F_k+n)\rfloor^s \
\displaybreak[2]\\
&= A(k,s,j)+F_k^j(F_{k+1}-\epsilon_k)^s+\sum_{n=1}^{F_{k-1}-1} (F_k+n)^j(F_{k+1}+\lfloor\phi n\rfloor)^s
\displaybreak[2]\\
&= A(k,s,j)+F_k^j\sum_{i=0}^s\binom siF_{k+1}^i(-\epsilon_k)^{s-i}
\\ &\qquad
+ \sum_{n=1}^{F_{k-1}-1} \sum_{\ell=0}^j\binom j\ell F_k^\ell n^{j-\ell}\sum_{i=0}^s\binom siF_{k+1}^i\lfloor\phi n\rfloor^{s-i}
\displaybreak[2]\\
&= A(k,s,j)+\sum_{i=0}^s\binom si(-\epsilon_k)^{s-i}F_k^jF_{k+1}^i
\\ &\qquad
+ \sum_{\ell=0}^j\sum_{i=0}^s\binom j\ell\binom siF_k^\ell F_{k+1}^iA(k-1,s-i,j-\ell).
\end{align*}
The above reduction, the identity $A(k,0,0)=F_k-1$, and the induction on $k+j+s$ implies that
$$
A(k,s,j)\in\operatorname{span}\{(\phi^i)^k,\,(-\phi^{i})^k:|i|\le j+s+1\};
$$
in particular, for a fixed choice of $s,j$ the sequence $\{A(k,s,j)\}_{k\ge1}$ is linearly recurrent of order at most $4(s+j)+6$.
Here we used this observation together with the following facts:
\begin{itemize}
\item If $\boldsymbol u=\{u_n\}_{n\ge 0}$ is a linearly recurrent sequence whose roots are all simple in some set $U$ then, for fixed integers $p$ and~$q$,
the sequence $\{u_{pn+q}\}_{n\ge 0}$ is linearly recurrent with simple roots in $\{\alpha^p: \alpha\in U\}$.
\item If $\boldsymbol u=\{u_n\}_{n\ge 0}$ and $\boldsymbol v=\{v_n\}_{n\ge 0}$ are linearly recurrent whose roots are all simple in some sets $U$ and $V$, respectively, then
$\boldsymbol{uv}=\{u_nv_n\}_{n\ge 0}$ is linearly recurrent whose roots are all simple in $UV=\{\alpha\beta: \alpha\in U,\; \beta\in V\}$.
\end{itemize}
In this context, the roots of a linearly recurrent sequence are defined as the zeroes of its characteristic polynomial, counted with their multiplicities.

It then follows that, for a fixed $s$, each of the sequences $\{A'(k,s)\}_{k\ge1}$ and $\{A'(k,s)\}_{k\ge1}$ are linearly recurrent of order at most $4s+6$.

\section{The proofs of the lemmas}
\label{sec4}

We first establish Lemma~\ref{conj2}.
By the argument in Section~\ref{sec3}, both $A(k,1)$ and $A'(k,1)$
are linearly recurrent with simple roots in the set $\{\pm \phi^l: |l|\le 2\}$.
The same is true for the right-hand sides in~\eqref{eq:conj2}. Since the set of roots is contained in a set with $10$ elements,
it follows that the validity of \eqref{eq:conj2} for $k=1,\ldots,10$ implies that the relations hold for all $k\ge 1$.

Lemmas \ref{conj3} and \ref{conj4} are similar. The argument in Section~\ref{sec3} shows that the left-hand sides
$\{A(k,3)\}_{k\ge 1}$ and $\{A'(k,3)\}_{k\ge 1}$ are linearly recurrent with simple roots in $\{\pm\phi^l: |l|\le 4\}$.
Splitting according to the parity of $k$ we deduce that
$\{A(2k,3)\}_{k\ge 1}$, $\{A(2k-1,3)\}_{k\ge 1}$, $\{A'(2k,3)\}_{k\ge 1}$ and $\{A'(2k-1,3)\}_{k\ge 1}$ are
linearly recurrent with simple roots in $\{\phi^{2l}: |l|\le 4\}$, a set with $9$~elements.
The same is true about the right-hand sides in the lemmas. Thus, if the relations hold for $k=1,\dots,9$ then they hold for all $k\ge 1$.

A few words about the computation. For the identities presented in Lemmas \ref{conj2}--\ref{conj4}, one can just use brute force to compute
$A(k,s)$ and $A'(k,s)$ for $s=1,3$ and $k=2\ell+i$, where $i\in \{0,1\}$ and $\ell=1,\ldots,9$, with any computer algebra system.
This takes a few minutes. If one would need to check it up to larger values of $k$, say around $100$,
the brute force strategy no longer works since the summation range up to $F_k-1$ becomes too large.
Instead one can use the recursion from Section~\ref{sec3} together with $A(k,0,0)=F_k-1$ to find subsequently $A(k,1,0)$, $A(k,2,0)$ and $A(k,3,0)$
for all desired $k$ and, similarly, $A(k,s,j)$ for small $j$ to evaluate $A'(k,s)$.

\section{The proof of Theorem \ref{conj1}}
\label{sec5}

Let us now address Theorem \ref{conj1}. When $k=4\ell$, this can be rewritten as
\begin{multline}
\label{eq:case4l}
F_{2\ell+1} ^2 L_{2\ell+2} L_{2\ell-1} (A'(4\ell,3) A(4\ell,1)^2  -  A(4\ell,3)A'(4\ell,1)^2)
\\
=A(4\ell,1)^2 A'(4\ell,1)^2(F_{2\ell+1}^2L_{2\ell+2} L_{2\ell-1}-1).
\end{multline}
Since $A(4\ell,s)$ and $A'(4\ell,s)$ are linearly recurrent (in $\ell$) with roots contained in $\{\phi^{4l}: |l|\le s+1\}$,
and both the left-most factor in the left-hand side and the right-most  factor in the right-hand side have each simple roots
in $\{\phi^{4l}: |l|\le 2\}$,  it follows that the both the left-hand side and the right-hand side are linearly recurrent
with simple roots contained in $\{\phi^{4l}: |l|\le 10\}$, a set with $21$ elements. Thus, if the above formula holds for $\ell=1,\dots,21$,
then it holds for all $\ell\ge 1$. A similar argument applies to the case when $k=4\ell+i$ for $i\in \{1,2,3\}$. Hence, all claimed formulas hold provided they hold for
all $k\le 100$, say.

Now we use the lemmas. For $k=4\ell$, Lemmas \ref{conj2}, \ref{conj3} and \ref{conj4} tell us that \eqref{eq:case4l},
after reducing the common factor $(F_{4\ell}-1)^2(F_{4\ell+1}-1)^2(F_{4\ell+2}-1)/16$, is equivalent to
\begin{align*}
&
F_{2\ell+1} ^2 L_{2\ell+2} L_{2\ell-1}
\times\Bigl(\frac15(F_{4\ell}-1)(L_{8\ell+4}-5L_{4\ell+3}+13)
\\ &\qquad\qquad\qquad\qquad\qquad
-(F_{4\ell+2}-1)(F_{4\ell-1}-1)(F_{4\ell+2}-1)\Bigr)
\\ &\quad
=(F_{4\ell}-1)^2(F_{2\ell+1}^2L_{2\ell+2} L_{2\ell-1}-1)
\end{align*}
(and one can perform further reduction using \eqref{eq:fact}).
It is sufficient to verify the resulting equality for $\ell=1,\dots,15$ and we have checked it for all $\ell=1,\dots,100$.
The remaining cases for $k$ modulo $4$ are similar.
We do not give further details here.

\section{The proof of Theorem~\ref{conj6}}
\label{sec6}

This follows from Lemma~\ref{conj2}, the classical formulas
\begin{equation}
\begin{alignedat}{2}
F_{4\ell}-1 &=  F_{2\ell+1} L_{2\ell-1}, &\qquad
F_{4\ell+1}-1 &= F_{2\ell} L_{2\ell+1}, \\
F_{4\ell+2}-1 &= F_{2\ell} L_{2\ell+2}, &\qquad
F_{4\ell+3}-1 &= F_{2\ell+2}L_{2\ell+1},
\end{alignedat}
\label{eq:fact}
\end{equation}
as well as known facts about the greatest common divisor of Fibonacci and Lucas numbers with close arguments. For example, for $k=2\ell$, we have
\begin{align*}
\operatorname{LCM}(2A(4\ell,1),2A'(4\ell,1))
&= \operatorname{LCM}((F_{4\ell+1}-1)(F_{4\ell}-1),(F_{4\ell+2}-1)(F_{4\ell}-1))
\\
&= \operatorname{LCM}(F_{2\ell} L_{2\ell+1}, F_{2\ell} L_{2\ell+2} ) F_{2\ell+1} L_{2\ell-1}
\\
&= F_{2\ell} L_{2\ell+1} L_{2\ell+2} F_{2\ell+1} L_{2\ell-1}
\\
&= F_{k+1} F_k L_{k+2} L_{k+1} L_{k-1},
\end{align*}
where we used the fact that $\gcd(L_{2\ell+1},L_{2\ell+2})=1$. The case $k=2\ell+1$ is similar.

\section{Further variations}
\label{sec7}

First, we give an informal account of a more general result lurking, perhaps, behind the formulas in Theorem~\ref{conj1}.
Consider a homogeneous (rational) function $r(\boldsymbol{x})=r(x_1,\dots,x_m)$ of degree~1, that is, satisfying
\begin{equation*}
r(t\boldsymbol{x}) = t  r(\boldsymbol{x}) \qquad\text{for}\quad t\in\mathbb Q,
\end{equation*}
and an algebraic number $\alpha$ solving the equation
\begin{equation}
\label{eq:charac}
\sum_{k=0}^m c_k \alpha^k=0,
\end{equation}
where $c_k$ are integers. If $r(\boldsymbol{x})$ vanishes at a vector $\boldsymbol{x}^*=(x_1^*,\dots,x_m^*)$ then we automatically have
\begin{equation}
\label{homog4}
\sum_{k=0}^m c_k r(\alpha^k\boldsymbol{x}^*)=0
\end{equation}
in view of the homogeneity of the function.
We can then inquire whether equation \eqref{homog4} is ``approximately'' true if $r(\boldsymbol{x}^*)=0$ is ``approximately'' true.
In this note, we merely examined the golden ration case in which \eqref{eq:charac} is $\alpha^2-\alpha-1=0$, while the choice
$$
r(x_1,\dots,x_m)=\frac{\sum_{n=1}^mx_n^3}{\bigl(\sum_{n=1}^mx_n\bigr)^2}
$$
for the rational function and $\boldsymbol{x}^*=(1,2,\dots,m)$ for its exact solution originated from the Nicomachus identity.

\medskip
Some further variations on the topic can be investigated in the $q$-direction, based on $q$-analogues of Nicomachus's theorem \eqref{eq:nico}
(consult with~\cite{Wa}).

\subsection*{Acknowledgements}
The first and third authors thank the Max Planck Institute for Mathematics (Bonn) for providing
excellent working conditions. The second author thanks his former PhD advisor Kenneth Stolarsky for communicating to this author
the conjectures that became a theorem in this note.

The first author was supported in part by NRF (South Africa) Grants CPRR160325161141 and an A-rated researcher award,
and by CGA (Czech Republic) Grant 17-02804S.
The second author was supported in part by the Institute of Mathematics of Universidad Autonoma de Santo Domingo,
through the grant FONDOCyT 2015-1D2-186, Ministerio de Educaci\'on Superior Ciencia y Tecnolog\'\i a (Dominican Republic).


\begin{thebibliography}{9}

\bibitem{Ki}
C. Kimberling,
``The Zeckendorf array equals the Wythoff array'',
\emph{Fibonacci Quart.} \textbf{33} (1995), 3--8.

\bibitem{St}
J. Stopple,
``A primer of analytic number theory. From Pythagoras to Riemann'',
Cambridge University Press, Cambridge (2003).

\bibitem{Wa}
S.\,O. Warnaar,
``On the $q$-analogue of the sum of cubes'',
\emph{Electronic J. Combin.} \textbf{11} (2004), no.~1, Note~13, 2~pp.

\end{thebibliography}
\end{document}